\documentclass{amsart}
\usepackage{array}
\usepackage{supertabular}

 \newcommand{\Pn}{\mathbb{P}^n}
 \newcommand{\PN}{\mathbb{P}^N}
 \renewcommand{\O}{\mathcal{O}}

 \DeclareMathOperator{\indet}{indet}

 \DeclareMathOperator{\rank}{rank}

 \newcommand{\A}{\mathbb{A}^1}

 \DeclareMathOperator{\Spec}{Spec}
 \DeclareMathOperator{\Sym}{Sym}

 \DeclareMathOperator{\Pic}{Pic}
 \DeclareMathOperator{\Hom}{Hom}

 \renewcommand{\P}{\mathbb{P}}

\begin{document}

 \newtheorem{thm}{Theorem}[section]
 \newtheorem{cor}[thm]{Corollary}
 \newtheorem{lem}[thm]{Lemma}
 \newtheorem{prop}[thm]{Proposition}

 \theoremstyle{definition}
 \newtheorem{defn}[thm]{Definition}
 \newtheorem{exam}[thm]{Example}

 \theoremstyle{remark}
 \newtheorem{obs}[thm]{Observation}
 \newtheorem{rem}[thm]{Remark}
 \newtheorem{fact}[thm]{Fact}
 \newtheorem{notation}[thm]{Notation}
 \numberwithin{equation}{section}


\fontsize{10}{12}\selectfont

\title[Morphisms of Hypersurfaces]
 {Towards Characterizing Morphisms\\ Between High Dimensional Hypersurfaces}

\author{ David C. Sheppard  }

\address{Department of Mathematics, MIT }

\email{sheppard@math.mit.edu}

\begin{abstract}
We prove the following theorem over an algebraically closed field of characteristic zero.
Let $f:X\rightarrow Y$ be a nonconstant morphism of hypersurfaces in $\mathbb P^n$, $n\geq 4$.  If $Y$ is nonsingular and of general type, then there is a morphism $F:\mathbb P^n\rightarrow\mathbb P^n$ such that $F|_X=f$ and $F^{-1}(Y)=X$.  As a corollary, we see that $\deg Y$ divides $\deg X$ with quotient $m$, and
$f$ is given by polynomials of degree $m$.
\end{abstract}
\maketitle

\section*{Survey of the Literature}

Classically, algebraic geometry has sought to classify varieties.  Recently focus has expanded to include the classification of morphisms between algebraic varieties.  Our goal in this paper is to shed some light on what type of morphisms can occur between hypersurfaces.
Let us give some results from the literature.

In \cite{PS}, Paranjape and Srinivas show that every nonconstant
morphism between smooth quadric hypersurfaces in $\Pn$ is an
isomorphism for $n\geq 4$.

Schuhmann shows in \cite{S1} that the degree of a
morphism from a smooth hypersurface of degree $d$ in $\P^4$ to a
smooth quadric threefold is bounded from above in terms of $d$.
For $d=3$ she obtains a very good bound on the degree of possible morphisms and proves that every morphism from a smooth cubic
threefold to a smooth quadric threefold is constant.

Amerik proves in \cite{A} that the degree of
a morphism from a smooth hypersurface in $\Pn$ to a smooth quadric
hypersurface in $\Pn$ is bounded from above for $n\geq 4$ using a
different argument than \cite{S1}.

Beauville shows in \cite{B} that every endomorphism of a smooth
hypersurface in $\Pn$ of dimension at least 2 and degree at least
3 is an automorphism.  For this he uses a Hurwitz-type inequality
from \cite{ARV}. We will generalize this inequality in Section 1.

\section*{Main Results and Point of View}

The main result of this paper is the following Theorem, which we prove in Section 2.

\subsection*{Theorem 2}\emph{
Assume the base field is algebraically closed of characteristic zero.  Let $f:X\rightarrow Y$ be a morphism of hypersurfaces in $\mathbb P^n$, $n\geq 4$.  If $Y$ is nonsingular and of general type, then there is a morphism $F:\mathbb P^n\rightarrow\mathbb P^n$ such that $f=F|_X$ and $F^{-1}(Y)=X$.
}
\vspace{.35cm}

The proof of Theorem 2, and our study of morphisms between hypersurfaces in general, relies on the following definition.

\subsection*{Definition}
 Let $f:X\rightarrow Y$ be a morphism of
projective $k$-varieties with specified very ample invertible
sheaves $\O_X(1)$ and $\O_Y(1)$ such that $f^*\O_Y(1)=\O_X(m)$.
Assume $X$ is positive dimensional so that there is a unique such
$m$.  We will refer to $m$ as the \emph{polynomial degree} of $f$,
because $f$ is given by polynomials of degree $m$.
\vspace{.25cm}

As a consequence of the Grothendieck-Lefschetz Theorem on the Picard group, every morphism $f:X\rightarrow Y$ between hypersurfaces of dimension at least 3 has a polynomial degree, cf. Lemma 2.1.  We will see that the conclusion of Theorem 2 is equivalent to the following statement:  $\deg Y$ divides $\deg X$ with quotient equal to the polynomial degree of $f$.

In light of this restatement, our first goal will be to bound the polynomial degree of possible morphisms between two given hypersurfaces. To do this, we generalize the Hurwitz-Type Inequality of \cite{ARV} in the case of morphisms between complete intersections.  This inequality gives good bounds on the polynomial degree of possible morphisms between hypersurfaces of high degree.  In particular, if $f:X\rightarrow Y$ is a morphism between hypersurfaces in $\P^n$ and $Y$ is nonsingular of degree at least $n+2$, i.e. Y is of general type, then the bound on the polynomial degree of $f$ is good enough to prove Theorem 2.

\subsection*{Acknowledgements} I would like to thank my advisor, A.
Johan de Jong, for his enumerable insights, gentle corrections,
and tireless enthusiasm during every stage of this project. Thanks
also to Roya Beheshti for helpful conversations about mapping
surfaces to threefolds, which lead to Proposition 2.5.

\section{A Hurwitz-Type Inequality}

In this section the ground field is algebraically closed of arbitrary characteristic.
The main result of this section is the following Theorem.

\subsection*{Theorem 1.}
\emph{Let $X$ be a complete intersection variety in $\Pn$ and $Y$ a
nonsingular projective variety of the same dimension as $X$.  Fix
a very ample invertible sheaf $\O_Y(1)$ on $Y$.  If
$f:X\rightarrow Y$ is a morphism such that
$f^*\O_Y(1)=\O_X(m)$ for some positive integer $m$ and the extension of function fields $k(Y)\rightarrow k(X)$ is separable, then
$$f^*c_\text{top}\left(\Omega^1_Y(2)\right)\leq c_\text{top}\left(\Omega^1_X(2m)\right).$$}

Theorem 1 is more general than the Hurwitz-Type Inequality
of \cite{ARV} in the sense that we do not assume the ground field has
characteristic zero or that $X$ is smooth.  However, it is less general in the sense that we assume $X$ is a complete intersection.

It is worth noting that Theorem 1 can fail if $Y$ is singular.  For example, $Y$ could be the image of a hyperplane in $\mathbb P^n$.

Also note that if $X$ is singular, then $\Omega^1_X(2m)$ is not locally free.  However, the top Chern class
$c_{top}(\Omega^1_X(2m))$ is defined via a finite locally free
resolution of $\Omega^1_X(2m)$, such as the conormal sequence for
$X\subset \Pn$.

We give some preliminary lemmas before proving Theorem 1.

\begin{lem}
If $X\subset \Pn$ is a positive dimensional complete intersection
and $f:X\rightarrow \PN$ is a morphism such that
$f^*\O_{\PN}(1)=\O_X(m)$ for some positive integer $m$, then $f$
extends to a rational map $F:\Pn\dashrightarrow\PN$ defined on a
Zariski open set containing $X$.
\end{lem}
\begin{proof}
Let $\xi_1,\dots,\xi_c$ be the homogeneous polynomials that
generate the homogeneous ideal of $X$, where $c$ is the
codimension of $X$ in $\Pn$.  Let $\xi_i$ have degree $d_i$, and
let $X_i:=V(\xi_1,\dots,\xi_i)$ so that $X=X_c\subset\dots\subset
X_0=\Pn$.  The morphism $f$ is given by sections $f_0,\dots,
f_N\in H^0(X_c, \O_{X_c}(m))$.  To lift $f$ from $X_c$ to a
rational map on $X_{c-1}$ we need to see that the restriction map
 $$
 H^0(X_{c-1},\O_{X_{c-1}}(m))\longrightarrow
 H^0(X_{c},\O_{X_{c}}(m))
 $$
is surjective.  So it suffices to check that
$H^1(X_{c-1},I(m))=0$, where $I\subset \O_{X_{c-1}}$ is the ideal
sheaf of $X_c$ in $X_{c-1}$. Since $I=\O_{X_{c-1}}(-d_c)$ is a
twisted structure sheaf and $X_{c-1}$ is a complete intersection
in $\P^n$ of dimension at least 2,
$$H^1(X_{c-1}, I(m))=0.$$ So the global sections $f_i$ lift from
$\O_{X_c}(m)$ to $\O_{X_{c-1}}(m)$. Continuing, we lift the $f_i$
to global sections $F_i$ of $\O_{\Pn}(m)$.  Set
$F=(F_0,\dots,F_N):\Pn\dashrightarrow\PN$, and note that $F$ is
undefined on $V(F_1,\dots,F_N)$, which is disjoint from $X$
because $F|_X=f$ is a morphism.
\end{proof}

The following positivity result essentially appears in \cite{F}.

\begin{lem}
Consider the following fiber square
$$
\begin{array}{ccc}
 W                  &
 \longrightarrow    &
 V                  \\
 \downarrow         &
                    &
 \downarrow         \\
 S                  &
 \stackrel{\rho}{\longrightarrow}    &
 T
\end{array}
$$
where $\rho:S\rightarrow T$ is a regular imbedding of varieties of
codimension $i$ and $V$ is a $k$-dimensional variety.  If
$N_ST$ is globally generated, then
$$\rho^*[V]=\sum_j \mu_j[R_j]+P$$
where the $R_j$ are the reduced structures on the
$(k-i)$-dimensional irreducible components of $W$, the $\mu_j$ are
positive integers, and $P$ is an effective $(k-i)$-cycle on $W$.
\end{lem}
\begin{proof}
We apply a positivity result to the Basic
Construction in \cite[Chapter 6]{F}.

Let $N$ be the pullback of $N_ST$ to $W$.  Then $N$ is a globally
generated vector bundle of rank $i$ on $W$, and the normal cone
$N_WV\hookrightarrow N$ is a purely $k$-dimensional closed
subscheme of $N$.  If $\sigma$ is the zero section of $N$, then
$\rho^*[V]$ is defined to be $\sigma^*[N_WV]$.  Recall that
$[N_WV]$ is the sum of the $k$-cycles associated to the
irreducible components of $N_WV$ taken with appropriate
multiplicities.

If $Z_j$ is a $(k-i)$-dimensional irreducible component of $W$,
then $N_WV$ has an irreducible component $C_j$ that dominates
$Z_j$.  If $R_j$ is the reduced structure on $Z_j$, then
$N|_{R_j}$ is the reduced structure on $C_j$ because $N|_{R_j}$ is
reduced and irreducible, $\dim C_j=k=\dim N|_{R_j}$, and $N|_{R_j}$
contains the reduced structure $(C_j)_{\text{red}}$.  So by
definition, $[C_j]=\mu_j[N|_{R_j}]$, where $\mu_j$ is the length
of the stalk of $\O_{C_j}$ at the generic point of $C_j$.
Therefore $\sigma^*[C_j]=\mu_j[R_j]$ because
$\sigma^*[N|_{R_j}]=[R_j]$. This accounts for the term $\sum
\mu_j[R_j]$ in the formula for $\rho^*[V]$.

Moreover, $\sigma^*$ takes effective cycles to effective cycles
because $N$ is globally generated, cf. \cite[Theorem 12.1(a)]{F}.
So if $C$ is an irreducible component of $N_WV$ other than one of
the $C_j$ described above, then $\sigma^*[C]$ is effective.  These other components $C$
account for $P$.
\end{proof}

\begin{lem}
Let $E$ be a globally generated vector bundle over a variety $X$.
Let $K\subseteq E$ be a closed subscheme with $\dim K < \rank E$.
Then there is a section $\sigma$ of $E$ such that $\sigma(X)\cap
K$ is empty.
\end{lem}
\begin{proof}
Since $E$ is globally generated, there is a surjective morphism of
vector bundles $\pi:X\times\mathbb{A}^{h^0}\rightarrow E$, where
$h^0:=h^0(X,E)$.  All the fibers of $\pi$ are affine spaces of the
same dimension.  So $\dim \pi^{-1}(K) = \dim K + (h^0-\rank E)$.
In other words, $h^0-\dim \pi^{-1}(K) = \rank(E)-\dim K>0$, whence
there is a constant global section $\tau$ of
$X\times\mathbb{A}^{h^0}$ over $X$ that does not intersect
$\pi^{-1}(K)$.  Take $\sigma = \pi\circ\tau$.
\end{proof}

\begin{lem}
Let $0\rightarrow L\rightarrow E\rightarrow F\rightarrow 0$ be a
short exact sequence of vector bundles on a complete variety $X$
such that $\rank F=\dim X$.  If $E$ is globally
generated, then for any morphism of vector bundles
$i:L\rightarrow E$ we have
$$\sigma^*[{i(L)}]\leq c_\text{top}(F)$$
where $\sigma$ is the zero section of $E$ and ${i(L)}$ is the
scheme-theoretic image of $i:L\rightarrow E$.  Equality holds if
$i$ is a closed immersion.
\end{lem}
\begin{proof}
First assume that $i$ is the given closed immersion $L\rightarrow
E$. Consider the fiber diagram
$$
\begin{array}{ccccc}
 & & L & \longrightarrow & X \\
 & & \downarrow\rlap{$i$} & & \downarrow\rlap{$\tau$} \\
 X & \stackrel{\sigma}{\longrightarrow} & E & \stackrel{\phi}{\longrightarrow} & F
\end{array}
$$
where $\tau$ is the zero section and $\phi$ is the quotient map.
Calculate
\begin{align*}
\sigma^*[i(L)]  &= \sigma^*\phi^*[X] \\
&= (\phi\circ\sigma)^*[X] \\
&=c_{top}(F).
\end{align*}
If $i:L\rightarrow E$ is any closed immersion of vector bundles on $X$ with quotient bundle
$F_i$, then the last statement of the
lemma follows from
$$c_{top}(F)=\Bigg\{\frac{c(E)}{c(L)}\Bigg\}_0=c_{top}(F_i).$$

Now assume that $i:L\rightarrow E$ is any morphism of vector bundles on $X$. If $\dim
{i(L)} <\dim L$, then $\sigma^*[{i(L)}] =0$. So
$\sigma^*[{i(L)}]\leq c_{top}(F)$ because $F$ is globally
generated, cf. \cite[Theorem 12.1(a)]{F}.  Therefore, we may assume
$\dim {i(L)}=\dim L$.

Let $\Sigma$ be the zero section of $E\times\A$ over $X\times\A$,
and let $\sigma_t:=\Sigma|_{X\times t}$ be the zero section of $E\times t$ over $X\times t$.  Let $x_t:X\times t\rightarrow X\times \A$ and $e_t:E\times
t\rightarrow E\times \A$ be the inclusion maps.  These maps fit together in the following fiber square.
\[
\begin{array}{ccc}
 X\times t & \stackrel{\sigma_t}\longrightarrow & E\times t \\
  \llap{$x_t$} \downarrow & & \downarrow \rlap{$e_t$} \\
 X\times\A & \stackrel\Sigma\longrightarrow & E\times \A
\end{array}
\]
Therefore, $x_t^* \Sigma^* \alpha = \sigma_t^*e_t^*\alpha$
for any cycle $\alpha$ on $E\times \A$.

Let $i_0=i$, and let $i_1$ be the given closed immersion $L\rightarrow E$ with quotient $F$.  Consider the morphism
\begin{align*}
I:L\times\mathbb A^1 &\longrightarrow E\times\mathbb A^1\\
(v,t) &\longrightarrow t\cdot i_1(v)+(1-t)\cdot i_0(v)
\end{align*}
of vector bundles on $X\times \mathbb A^1$, and let $i_t=I|_{L\times t}:L\times t\rightarrow E\times t$.
Let $Z$ denote the scheme-theoretic image of $I$, and define
$$
\lambda_t := x_t^*\Sigma^*[Z] = \sigma^*_te_t^*[Z].
$$
Since  $i_t:L\rightarrow E$ is a closed immersion for
all $t$ in some neighborhood of $1\in \A$,  $i_1(L)=e_1^{-1} (Z)$.  So $e_1^{*}[Z]=[i_1(L)]$, which implies
 $\lambda_1 = \sigma_1^*[i_1(L)]$.

Note that ${i_0(L)}$ is an irreducible component of
$e_0^{-1}(Z)$ because $i_0(L)\subset e_0^{-1}(Z)$ and
$$\dim e_0^{-1}(Z)+1 = \dim Z = \dim {i_0(L)}+1.$$
So by Lemma 1.2, $e_0^*[Z] =[{i_0(L)}]+P$ for some effective
cycle $P$. Therefore,
$$
\lambda_0=\sigma_0^*[{i_0(L)}]+\sigma_0^*P
$$
and $\sigma_0^*P$ is an effective 0-cycle because $E$ is globally
generated.

If $\alpha$ is any cycle on $X\times \A$, then since $X\times\A\rightarrow\A$ is proper, the degree of the restriction
$x_t^*\alpha$ of $\alpha$ to the fiber $X\times t$ does not depend
on $t$, cf.
\cite[Proposition 10.2]{F}.
Take $\alpha=\Sigma^*[Z]$ to see that the 0-cycles $\lambda_t$ on $X$ all have the same degree.  Now we can calculate
\begin{align*}
\deg \sigma_0^*[{i_0(L)}] & \leq \deg \lambda_0 \\
 & = \deg \lambda_1 \\
 & = \deg \sigma_1^*[i_1(L)] \\
 & = c_{top}(F).
\end{align*}
\end{proof}

\subsection*{Proof of Theorem 1.}
Let $f:X\rightarrow Y$ be a morphism of projective varieties.
Assume $X$ is a complete intersection in $\Pn$ and $Y$ is
nonsingular of the same dimension as $X$.  Fix a very ample
invertible sheaf $\O_Y(1)$ on $Y$ with corresponding projective
embedding $Y\hookrightarrow \PN$.  Assume $f^*\O_Y(1)=\O_X(m)$
for some positive integer $m$, which implies that $f$ is finite
and surjective. Assume also that $f$ is a separable morphism.

By Lemma 1.1, there is a rational map $F:\Pn\dashrightarrow\PN$
defined on a Zariski open subset of $\Pn$ containing $X$ such that
$f=F|_X$.  This extended map $F$ induces a morphism
$f^*(\Omega^1_{\PN}|_Y)\rightarrow\Omega^1_{\Pn}|_X$, which gives
the following commutative diagram of sheaves on $X$:

\begin{equation}
\begin{array}{ccccccccc}
 0                      &
 \longrightarrow        &
 f^*(I_Y/I_Y^2(2))      &
 \longrightarrow        &
 f^*(\Omega^1_{\PN}|_Y (2))  &
 \longrightarrow        &
 f^*(\Omega^1_Y (2))         &
 \longrightarrow        &
 0                      \\
                        &
                        &
 \downarrow             &
                        &
 \downarrow             &
                        &
 \downarrow             &
                        &
                        \\
 0                      &
 \longrightarrow        &
 I_X/I_X^2(2m)          &
 \longrightarrow        &
 \Omega^1_{\Pn}|_X (2m) &
 \longrightarrow        &
 \Omega^1_X (2m)        &
 \longrightarrow        &
 0
\end{array}
\end{equation}
\smallskip
where $I_X\subset\O_{\Pn}$ and $I_Y\subset\O_{\PN}$ are the ideal
sheaves of $X$ and $Y$.  The bottom row is exact because $I_X$ is
the ideal sheaf of a reduced complete intersection.

To apply our intersection-theoretic lemmas, we transform diagram
(1.1) of sheaves on $X$ into a diagram of schemes over $X$ by
applying the covariant functor
\begin{align*}
\Phi : \{\text{coherent sheaves on} X\} &\longrightarrow\{\text{schemes of finite type over} X\}\\
\mathcal{F} &\longrightarrow \Spec
\left(\Sym_{\O_X}[\Hom_{\O_X}(\mathcal{F},\O_X)]\right)
\end{align*}
where $\Sym_{\O_X}(-)$ denotes the symmetric algebra of an
$\O_X$-module.
If $\mathcal F$ is a locally free sheaf, then $\Phi(\mathcal F)$ is the vector bundle whose sheaf of sections is $\mathcal F$.  Apply $\Phi$ to diagram (1.1), and
denote the resulting diagram of $X$-schemes by
\begin{equation}
\begin{array}{ccccc}
 L_Y                                    &
 \stackrel{i_Y}{\longrightarrow}        &
 E_Y                                    &
 \rightarrow                            &
 F_Y                                    \\
 \downarrow                             &
                                        &
 \downarrow \rlap{$\psi$}               &
                                        &
 \downarrow                             \\
 L_X                                    &
 \stackrel{i_X}{\longrightarrow}        &
 E_X                                    &
 \rightarrow                            &
 F_X
\end{array}
\end{equation}
Note that every scheme in (1.2) is a vector bundle on $X$,
except for $F_X$ if $X$ is singular. Also note that $i_Y$ is a
closed immersion because $Y$ is nonsingular, and that $E_X$, $E_Y$
are generated by global sections because $\Omega^1_{\Pn}(a)$ is globally generated for $a\geq 2$.

Let $\sigma$ be the zero section of $E_Y$ so that
$\psi\circ\sigma$ is the zero section of $E_X$.  By Lemma 1.4,
\begin{equation}
f^*c_{top}\left(\Omega^1_Y(2)\right)=\sigma^*[i_Y(L_Y)].
\end{equation}
Let $i_X(L_X)$ be the scheme-theoretic image of
$i_X:L_X\rightarrow E_X$.  By equation (1.3), it suffices to show
\begin{align}
\sigma^*[i_Y(L_Y)] &\leq (\psi\circ\sigma)^*[i_X(L_X)] \\
 & \leq c_{top}\left(\Omega^1_X(2m)\right).
\end{align}

To prove (1.4) it is enough to show
$\psi^*[i_X(L_X)]=[i_Y(L_Y)]+P_Y$ for some effective cycle $P_Y$
on $E_Y$.  Indeed, $\sigma^*P_Y$ is effective because $E_Y$ is globally generated, whence
\begin{align*}
\sigma^*[i_Y(L_Y)] &\leq \sigma^*[i_Y(L_Y)] +\sigma^*P_Y \\
&=\sigma^*\psi^*[i_X(L_X)] \\
& = (\psi\circ\sigma)^*[i_X(L_X)].
\end{align*}
Consider the fiber diagram
$$
\begin{array}{ccccc}
 \psi^{-1}(i_X(L_X))     &
 \longrightarrow                    &
 E_Y\times_Xi_X(L_X)     &
 \longrightarrow                    &
 i_X(L_X)                \\
 \downarrow                         &
                                    &
 \downarrow                         &
                                    &
 \downarrow                         \\
 E_Y                                &
 \stackrel{\Gamma_\psi}{\longrightarrow}        &
 E_Y\times_X E_X                    &
 \stackrel{\pi_2}{\longrightarrow}               &
 E_X
\end{array}
$$
where $\Gamma_\psi$ is the graph of $\psi$.  Then $\psi =
\pi_2\circ\Gamma_\psi$, whence
\begin{align*}
\psi^*[i_X(L_X)] &= \Gamma_\psi^*\pi_2^*[i_X(L_X)] \\
& =\Gamma_\psi^*[E_Y\times_Xi_X(L_X)].
\end{align*}
Since $\Gamma_\psi$ is a section of the vector bundle $E_Y\times_X
E_X $ over $E_Y$
$$
\Gamma_\psi^*(N)\cong E_Y\times_XE_X
$$
where $N$ is the normal bundle of $\Gamma_\psi(E_Y)$ in $ E_Y\times_XE_X$.  Therefore, $N$ is globally generated because $E_X$ is globally generated over $X$.  So it suffices by Lemma
1.2 to show that $i_Y(L_Y)$ is an irreducible component of
$\psi^{-1}(i_X(L_X))$.

By the assumption that $k(Y)\hookrightarrow k(X)$ is a separable
field extension, the stalk of $\Omega_{X/Y}^1$ at the generic
point of $X$ is $\Omega^1_{k(X)/k(Y)}=0$.  Hence there is some
nonempty open $U$ in $X$ such that the restriction of
$f^*\Omega_Y^1\rightarrow\Omega_X^1$ to $U$ is an isomorphism of
locally free sheaves. So when diagram (1.1) is restricted to $U$, the morphism
$f^*\left(\Omega_Y^1(2)\right)\rightarrow\Omega_X^1(2m)$ becomes an isomorphism.  Hence $F_Y\rightarrow F_X$ is an isomorphism when restricted to $U$.  It follows that $\psi^{-1}(i_X(L_X))$
and $i_Y(L_Y)$ coincide over $U$. Therefore, $i_Y(L_Y)$ is an irreducible
component of $\psi^{-1}(i_X(L_X))$. This establishes equation
(1.4).

\smallskip
To prove equation (1.5), it suffices by Lemma 1.4 to show that
there is a closed immersion $L_X\rightarrow E_X$ of vector
bundles, i.e. that there is a morphism of locally free sheaves
$I_X/I_X^2\rightarrow \Omega^1_{\Pn}|_X$ with empty degeneracy
locus.

Let $X$ be cut out by homogeneous polynomials $\xi_1,\dots,\xi_c$
where $c$ is the codimension of $X$ in $\Pn$.  Let $a_i=\deg
\xi_i$ so that
$$I_X/I_X^2 \cong \bigoplus_{i=1}^c\O_X(-a_i).$$
By decreasing $n$ if necessary, we assume $a_i\geq 2$ for
each $a_i$.   If $c=0$,
then there is nothing to prove, so assume $c>0$. We will construct
a morphism $\bigoplus\O_X(-a_i)\rightarrow \Omega^1_{\Pn}|_X$ with
empty degeneracy locus one summand at a time.

Since $a_1\geq 2$, the locally free sheaf $\Omega^1_{\Pn}|_X(a_1)$
is globally generated.   By Lemma 1.3, the rank $n$ vector bundle
$\Phi\left(\Omega_{\P^n}^1|_X(a_1)\right)$ has a section that
avoids the zero section. Hence there is a
morphism $\sigma_1 : \O_X\rightarrow \Omega^1_{\Pn}|_X(a_1)$ with
empty degeneracy locus.  Tensoring $\sigma_1$ with $\O_X(-a_1)$
gives a morphism $\phi_1:\O_X(-a_1)\rightarrow \Omega^1_{\Pn}|_X$
with empty degeneracy locus.  If $c=1$, we are done.

If $c\geq 2$, then let $\phi_1^\prime:\O_X(a_2-a_1)\rightarrow \Omega^1_{\Pn}|_X(a_2)$ denote the morphism obtained from $\phi_1$ by tensoring with $\O_X(a_2)$. Since
 $n> n-c+1$ and the image of $\Phi(\phi_1^\prime)$ has dimension $n-c+1$, Lemma 1.3 implies that there is a section of
$\Phi\left(\Omega_{\P^n}^1|_X(a_2)\right)$ that avoids the image of
$\Phi(\phi_1^\prime)$.  In other words,
there is a morphism $\sigma_2:\O_X\rightarrow
\Omega_{\P^n}^1|_X(a_2)$ such that
$$
 \phi_1^\prime\oplus\sigma_2:\O_X(a_2-a_1)\oplus\O_X
 \rightarrow \Omega^1_{\P^n}|_X(a_2)
$$
has empty degeneracy locus. If $\phi_2:\O_X(-a_2)\rightarrow\Omega^1_{\mathbb P^n}|_X$ is obtained from $\sigma_2$ by tensoring with $\O_X(-a_2)$, then tensoring the above morphism with $\O_X(-a_2)$ yields a morphism
 $$
 \phi_1\oplus\phi_2:\O_X(-a_1)\oplus\O_X(-a_2)
 \longrightarrow \Omega^1_{\P^n}|_X
 $$
with empty degeneracy locus.

Continuing like this, we obtain a morphism
$\bigoplus \phi_i:\bigoplus\O_X(-a_i)\rightarrow \Omega^1_{\Pn}|_X$ with empty
degeneracy locus.  This completes the proof of Theorem 1.

\section{Morphisms Between Hypersurfaces}

We will apply Theorem 1 to the case of hypersurfaces in $\P^n$.  We fix the notation and assumptions of the following discussion for the rest of the paper.

Let $f:X_d\rightarrow Y_e$ be a nonconstant morphism of
hypersurfaces of the indicated degrees in $\Pn$, $n\geq 4$.  Assume $X$ is integral and $Y$ is nonsingular.  We also assume  $e\geq 3$ because the inequality of Chern classes in Theorem 1 only gives good information in this range.

The Grothendieck-Lefschetz Theorem, \cite[Theorem 4.3.2]{H1}, states that $\Pic X$ is generated by $\O_X(1)$.  Therefore, $f^*\O_Y(1)=\O_X(m)$ for some nonnegative integer $m$.  Since $f$ is not constant, the polynomial degree $m$ of $f$ is positive.

As $f^*\O(1)$ is ample, $f$ is finite. Therefore, $f$ induces a finite extension of function fields.  We assume that the extension $k(Y)\rightarrow k(X)$ is separable.

Now we introduce a hypersurface $H$ that will be central to our study of $f:X\rightarrow Y$.
By Lemma 1.1, the morphism $f:X\rightarrow Y$ of Theorem 2 extends to a rational
map $F:\P^n\dashrightarrow\P^n$ defined at all but finitely many
points away from $X$. Since $e\geq 3$, $Y$ is not the image of a hyperplane in $\mathbb P^n$, because the only smooth variety that is the image of a
morphism from a projective space is projective space itself, cf.
\cite{L}.  Therefore, $Y$ is not the image of $F:\mathbb P^n\dashrightarrow \P^n$.  Hence,  $F$ is dominant because its image is irreducible and contains $Y$.
It follows that $F^{-1}(Y)$ is a hypersurface in $\P^n$.  Since $X\subset F^{-1}(Y)$, we may define the hypersurface $H$ in $\mathbb P^n$ as the difference of divisors
$$H:= F^{-1}(Y)-X.$$
We will study $H$ because $F^{-1}(Y)=X$ if and only if  $H=0$ as a divisor on $\P^n$, i.e. $H$ is empty.

\subsection{First Calculations}
The ground field will be algebraically closed of arbitrary characteristic in this subsection.
We will derive closed formulas for
$c_{n-1}\left(\Omega_X^1(2m)\right)$ and
$f^*c_{n-1}\left(\Omega_Y^1(2)\right)$.
So consider the short exact sequences
$$
0\longrightarrow \O_X(-d)\longrightarrow \Omega_{\P^n}^1
\longrightarrow \Omega_X^1\longrightarrow 0
$$
$$
0\longrightarrow \Omega_{\P^n}^1 \longrightarrow
\O_{\P^n}(-1)^{\oplus n+1} \longrightarrow \O_{\P^n}
\longrightarrow 0
$$
Let $h:=c_1(\O_X(1))$, and calculate the total Chern class of $\Omega_X^1$ to be
\begin{align*}
c\left(\Omega_X^1\right) &= \frac{(1-h)^{n+1}}{1-dh} \\
 &= \left( \sum_{i=0}^{n-1} \binom{n+1}{i}(-h)^i\right) \cdot
    \left(\sum_{j=0}^{n-1} (dh)^j\right)
\end{align*}
The $i^{\text{th}}$ Chern class of $\Omega_X^1$ is therefore given by
$$
c_i\left(\Omega_X^1\right) =  h^i \cdot \sum _{j=0}^i
(-1)^j\binom{n+1}{j}d^{i-j}
$$
The usual calculation with Chern roots shows
\begin{align}
c_{n-1}\left(\Omega_X^1(2m)\right)
 &= \sum _{i=0}^{n-1} c_i(\Omega_X^1)(2mh)^{n-1-i} \\
 &= h^{n-1} \sum _{i=0}^{n-1} \sum _{j=0}^{i} (-1)^j\binom{n+1}{j}
 d^{i-j} (2m)^{n-1-i}
\end{align}

Notice that for each pair of integers $a,b$ such that $a\geq 0$, $b\geq 0$, and $a+b\leq n-1$, the monomial $d^{a+1}(2m)^b$ has coefficient $(-1)^N\binom{n+1}{N}$ in (2.2), where $N=n-1-a-b$.  Therefore, we introduce the notation
$$
\Phi_N(x,y):= x^N+x^{N-1}y+\dots+x y^{N-1}+y^N
$$
and use the observation $h^{n-1}=d$ to obtain
$$
c_{n-1}\left(\Omega_X^1(2m)\right)
= d \sum_{k=0}^{n-1}(-1)^k \binom{n+1}k \Phi_{n-1-k}(d, 2m).
$$
We continue the calculation of $c_{n-1}(\Omega_X^1(2m))$ as follows.
\begin{align*}
 c_{n-1}\left(\Omega_X^1(2m)\right) &=  \sum_{k=0}^{n}(-1)^k \binom{n+1}k d\left(\frac{d^{n-k}-(2m)^{n-k}}{d-2m}
  \right) \\
 &= \frac1{d-2m}\ \Bigg\{ \sum_{i=0}^{n}(-1)^i \binom{n+1}{i} d^{n+1-i} \\
  & \hspace{3 cm} - d \sum_{j=0}^{n}(-1)^j \binom{n+1}{j}
 (2m)^{n-j} \Bigg\} \\
 &= \frac1{(2m)(d-2m)}\ \Bigg \{ 2m \sum_{i=0}^{n}(-1)^i \binom{n+1}{i} d^{n+1-i}\\
   &\hspace{3 cm} - d \sum_{j=0}^{n}(-1)^j \binom{n+1}{j} (2m)^{n+1-j} \Bigg \}
    \\
 &= \frac1{(2m)(d-2m)}\  \Bigg \{
 2m \sum_{i=0}^{n+1}(-1)^i \binom{n+1}{i} d^{n+1-i} \\
    & \hspace{3 cm} - d \sum_{j=0}^{n+1} (-1)^j \binom{n+1}{j}
    (2m)^{n+1-j} \\
    & \hspace{3.3 cm} +(-1)^{n+1}(d-2m) \Bigg \} \\
 &= \frac{2m(d-1)^{n+1}-d(2m-1)^{n+1}+(-1)^{n+1}(d-2m)}
    {2m(d-2m)}
\end{align*}
Introducing $x=2m-1$ and $y=d-1$, we calculate
$c_{n-1}\left(\Omega_X^1(2m)\right)$ to be
\begin{align*}
c_{n-1}\left(\Omega_X^1(2m)\right)
 &= \frac{(x+1)y^{n+1}-(y+1)x^{n+1}+(-1)^{n+1}(y-x)}{(x+1)(y-x)} \\
 &= \frac{xy(y^n-x^n)+(y^{n+1}-x^{n+1})
       +(-1)^{n+1}(y-x)}{(x+1)(y-x)} \\
 &= \frac{xy\Phi_{n-1}(x,y)+\Phi_n(x,y) +(-1)^{n+1}}{x+1} \\
 &= \frac{xy\Phi_{n-1}(x,y) + x\left( \Phi_{n-1}(x,y)+ \frac{y^n}{x}\right)
 +(-1)^{n+1}}{x+1} \\
 &=
 \frac{x(y+1) \Phi_{n-1}(x,y)+y^n
 +(-1)^{n+1}}{x+1}.
\end{align*}
Therefore we obtain the formula
\begin{equation}
 \boxed{
  \hspace{0.1 cm}
  c_{n-1}\left(\Omega_X^1(2m)\right)
  =
  \frac{d(2m-1)\Phi_{n-1}(2m-1,d-1)+(d-1)^n+(-1)^{n+1}}{2m}
  \hspace{0.1 cm}
 }
\end{equation}
By taking $m=1$ and substituting $e$ for $d$ in formula (2.3), we
have a formula for $c_{n-1}\left(\Omega_Y^1(2)\right)$.  Therefore, we can use the equations
$$
f^*c^{n-1}\left(\Omega_Y^1(2)\right)=\deg f\cdot c^{n-1}\left(\Omega_Y^1(2)\right)
\hspace{.6cm} \text{and} \hspace{.6cm}
\deg f=\frac{dm^{n-1}}{e}
$$
to derive the following formula for $ f^*c_{n-1}\left(\Omega_Y^1(2)\right)$
\begin{equation}
 \boxed{
  \hspace{0.1 cm}
  f^*c_{n-1}\left(\Omega_Y^1(2)\right)
  =
  \frac{dm^{n-1}}{e} \left(
  \frac{e\ \Phi_{n-1}(1,e-1)+(e-1)^n+(-1)^{n+1}}{2}
  \right)
  \hspace{0.1 cm}
 }
\end{equation}

We will need the following polynomial fact in the proof of Proposition 2.2.

\begin{lem}
If $x$, $y$ are positive real numbers with $x\geq 3$, $N\geq
3$ is an integer, and $\Phi_N(y,2)>(x+1)^N+1$, then $y>x$.
\end{lem}
\begin{proof}
Since $\Phi_N(y,2)$ increases with respect to $y$ it suffices to
show that if $x\geq 3$, then $\Phi_N(x,2) \leq (x+1)^N+1$.  Notice
that the coefficients of the polynomial
$P(x)=(x+1)^N+1-\Phi_N(x,2)$ have only one sign change.  So by Descarte's rule of signs, $P(x)$ has only one positive real root.  Therefore,
since $P(0)<0$, it suffices to check that $P(3)\geq 0$.  One
easily checks this for $N\geq 3$.
\end{proof}

\begin{prop} In the notation established at the beginning of this Section:\\
(1) For each triple $(d,e,n)$ there is an integer $M=M(d,e,n)$ such that $m\leq M$. \\
(2) $d\geq e$.\\
(3) If $d=e$, then $m=1$.\\
(4) If $e\geq 5$, then $d-1> m(e-2)$.
\end{prop}

\subsection*{Remarks}
If the base field is $\mathbb C$ and $X$ is nonsingular, then (2) has the following proof, which is independent of Theorem 1.
There is an injection of singular cohomology rings $H^*(Y,\mathbb
C )\rightarrow H^*(X,\mathbb C)$. So in this case, (2) can be
proved by computing the dimension of the middle cohomology groups
of $X$ and $Y$.

Part (3) is a generalization of the result in \cite{B}
that in characteristic zero every endomorphism of a smooth
hypersurface of degree at least 3 and dimension at least 2 is an
automorphism. We only assume $X$ and $Y$ have the same degree, not
that $X=Y$, and we do
not assume characteristic zero, only that the morphism is separable. The case $n=3$ can be checked without much work using Theorem 1.

Part (4) will be needed for the proof of Theorem 2.
\begin{proof}
Theorem 1 states (2.3) $\geq$ (2.4).  Dividing both sides
of this inequality by $dm^{n-1}$ results in
 $$
 \frac{2m-1}{2m}\ \Phi_{n-1}\left(\frac{d-1}{m},\frac{2m-1}{m}\right)
 \ +\ \frac{d-1}{d}\ \frac{1}{2m}\left(\frac{d-1}{m}\right)^{n-1}+\ \frac{(-1)^{n+1}}{2m^nd}
 $$
 $$
 \geq\frac12\Phi_{n-1}(e-1, 1)+\frac{(e-1)^n+(-1)^{n+1}}{2e}
 $$
Using $\frac{d-1}{d}<1$, combine the first two terms in the above
inequality to see
\begin{align*}
 \Phi_{n-1}\left(\frac{d-1}{m}, 2\right) \
 & >\ \frac12\Phi_{n-1}(e-1, 1)+\frac{(e-1)^n+(-1)^{n+1}}{2e} \\
 &=\ \frac12 \left(\frac{(e-1)^n -1}{(e-1)-1}\right)
 +\frac{(e-1)^n+(-1)^{n+1}}{2e}\\
 & =\ \frac{2(e-1)^{n+1}-e+(-1)^{n+1}(e-2)}{2(e-2)} \\
 &\geq\ \frac{(e-1)^{n+1}-(e-1)}{e(e-2)}.
\end{align*}
Since we assume $e\geq 3$, this implies
\begin{equation}
\Phi_{n-1}\left(\frac{d-1}{m}, 2\right) > (e-1)^{n-1}+1.
\end{equation}

Suppose $m$ were not bounded from above.  Taking the limit of
(2.5) as $m\rightarrow\infty$ shows $2^{n-1}\geq (e-1)^{n-1}+1$.
This contradiction proves (1).

To prove (4), notice that if $e\geq 5$, then Lemma 2.2 and
inequality (2.5) imply that $\frac{d-1}{m}>e-2$.

If $m=1$, then $d=e$, as follows.  Use Lemma 1.1 to extend $f$
to a rational map $F:\P^n\dashrightarrow\P^n$ with
$F^*\O(1)=\O(1)$. The image of $F$ is a linear subspace of $\P^n$
that contains $Y$, namely $\P^n$ itself. So $F$ is an automorphism
of $\P^n$, and $d=e$.

To prove (2) and (3) we assume $m\geq 2$ and prove $d>e$.  If
$e\geq 5$, then $d>e$ by (4). The cases $e=3$ and $e=4$ can be
checked by hand  in case $n=4$, and it suffices to check (2)
and (3) for the case $n=4$ because the upper bounds on $m$ given
by the inequality of Theorem 1 improve as $n$ increases.
\end{proof}

\begin{cor}
Let $f:X\rightarrow Y$ be a nonconstant separable morphism of hypersurfaces  in $\P^n$, $n\geq 4$, such that $Y$ is nonsingular and $\deg X=\deg Y\geq 3$. There is an automorphism $F:\P^n\rightarrow\P^n$ such that $f=F|_X$.
\end{cor}
\begin{proof}
By Proposition 2.2(3), $f^*\O_Y(1)=\O_X(1)$.  By Lemma 1.1, there is a rational map $F:\P^n\dashrightarrow\P^n$ such that $f=F|_X$.  Since $F^*\O(1)=\O(1)$, the image of $F$ is a linear subspace of $\P^n$ containing $Y$.  So $F$ is in fact an automorphism.
\end{proof}

\subsection{Hypersurfaces of General Type}
We now assume the ground field is algebraically closed of characteristic zero.
The purpose of this subsection is to prove the following Theorem.

\vspace{.3cm}
\noindent\textbf{Theorem 2.}\emph{
Assume the base field is algebraically closed of characteristic zero.  Let $f:X\rightarrow Y$ be a morphism of hypersurfaces in $\mathbb P^n$, $n\geq 4$.  If $Y$ is nonsingular and of general type, then there is a morphism $F:\mathbb P^n\rightarrow\mathbb P^n$ such that $f=F|_X$ and $F^{-1}(Y)=X$.
}
\vspace{.35cm}

The proof will rely on Proposition 2.4 below, which is an inequality that will bound the polynomial degree $m$ of $f$ from below.

To prove Theorem 2, we focus our attention on the hypersurface $H:=F^{-1}(Y)-X$ in $\P^n$ defined at the beginning of Section 2.  In particular, we wish to show that $H$ is the $0$ divisor, i.e. that $F^{-1}(Y)=X$.

Define $\Sigma$ to be an irreducible component of a general hyperplane section of $H$, taken with its reduced structure.  Then $F$ is defined at every point of $\Sigma$ because $F$ is undefined at only finitely many points in $\P^n$. We will analyze the morphism $F|_\Sigma:\Sigma\rightarrow Y$ using the following Proposition.

\begin{prop}
Let $\Sigma$ be an integral hypersurface in $\P^{n-1}$ and $Y$ be
a smooth hypersurface in $\Pn$, $n\geq 4$.  Let $\delta=\deg \Sigma$, and $e=\deg
Y$.   If $g:\Sigma\rightarrow Y$ is a
morphism with $g^*\O_Y(1)=\O_\Sigma(m)$ for some positive integer
$m$, then
$$n-\delta+m(e-n) \leq 0.$$
\begin{proof}
First we claim that there is a canonical morphism
\begin{align}
\bigwedge^{n-2}\Omega^1_\Sigma\rightarrow\omega^o_\Sigma
\end{align}
that is an isomorphism on the nonsingular locus of $\Sigma$, where
$\omega^o_\Sigma$ is the dualizing sheaf of $\Sigma$.
Let $I$ denote the ideal sheaf of $\Sigma$ in $\Pn$.  Since
$\Sigma$ is a reduced local complete intersection in $\Pn$, there
is a short exact sequence
\begin{equation}
0\longrightarrow I/I^2\longrightarrow\Omega^1_{\Pn}|_\Sigma
\longrightarrow\Omega^1_\Sigma\longrightarrow 0.
\end{equation}
Therefore the morphism
 $$
 \Phi\ :
 \bigg(\bigwedge^2 I/I^2\bigg)\otimes\bigwedge^{n-2}\Omega^1_\Sigma
 \longrightarrow
 \bigwedge^n\Omega^1_{\Pn}|_\Sigma
 $$
 $$
 \xi_1\wedge\xi_2\otimes
 d\overline{\phi}_1\wedge\dots\wedge d\overline{\phi}_{n-2}
   \longmapsto
  d\xi_1\wedge d\xi_2\wedge
 d\phi_1\wedge\dots\wedge d\phi_{n-2}
 $$
is well-defined.  Using the formula
 $$
  \omega^o_\Sigma \cong \bigg(\bigwedge^2I/I^2\bigg)^{-1}\otimes
  \bigwedge^n\Omega^1_{\Pn}|_\Sigma
 $$
tensor $\Phi$ with the dual of the invertible sheaf $\bigwedge^2
I/I^2$ to obtain the morphism (2.6).  This
is an isomorphism when restricted to
$\Sigma_{reg}$, because all the sheaves in (2.7) are locally free on $\Sigma_{reg}$.

Since $g^*\O_Y(1)=\O_\Sigma(m)$ is ample, $g$ has finite fibers.
So the canonical morphism
$g^*\Omega^1_Y\rightarrow\Omega^1_\Sigma$ is a surjection at the
generic point of $\Sigma$ by the characteristic zero assumption of
this subsection. By taking exterior powers and composing with (2.6), we obtain a composite morphism
$$
\bigwedge^{n-2}g^*\Omega^1_Y\longrightarrow\bigwedge^{n-2}\Omega^1_\Sigma
\longrightarrow\omega_\Sigma^o
$$
that is a surjection at the generic point of $\Sigma$. Since
$\Sigma$ is a hypersurface in $\mathbb P^{n-1}$ of degree $\delta$,
$\omega_\Sigma^o=\O_\Sigma(\delta-n)$.  So dualizing the above
morphism gives the exact sequence \vspace{-0.2 cm}
 \begin{equation}
 0 \longrightarrow \O_\Sigma(n-\delta)\longrightarrow\bigwedge^{n-2}g^*T_Y.
 \end{equation}
This is an injection because it is an injection at the
generic point of $\Sigma$ and $\O_\Sigma(n-\delta)$ is
torsion-free.  Tensoring (2.8) with $\O_\Sigma(m(e-n)-1)$ and applying the
formula $\bigwedge^{n-2}T_Y = \Omega_Y^1(-K_Y)$ yields the exact
sequence
\begin{equation}
 0\longrightarrow \O_\Sigma(n-\delta +m(e-n)-1) \longrightarrow
 \left(g^*\Omega^1_Y \right)(m-1).
\end{equation}

Tensoring with $\O_Y(-K_Y)=\O_Y(n+1-e)$, the conormal sequence for
$Y$ in $\P^n$ and the Euler sequence for $\P^n$ give
the following short exact sequences,  respectively:
\begin{equation}
 0\longrightarrow \O_Y(n+1-2e) \longrightarrow
 \Omega^1_{\P^n}\otimes\O_Y(-K_Y) \longrightarrow \Omega_Y^1(-K_Y) \longrightarrow0
\end{equation}
\begin{equation}
 0\longrightarrow \Omega^1_{\P^n}\otimes\O_Y(-K_Y) \longrightarrow
 \O_Y(n-e)^{\oplus n+1} \longrightarrow\O_Y(n+1-e)\longrightarrow 0
\end{equation}

Tensor (2.10) and (2.11) with $\O_Y(e-n)$, apply $g^*$, then
tensor with $\O_\Sigma(-1)$ to obtain the following short exact
sequences:
\begin{equation}
 0\longrightarrow \O_\Sigma(m(1-e)-1) \longrightarrow
 \left(g^*\Omega^1_{\P^n}\right)(m-1) \longrightarrow \left(g^*\Omega_Y^1\right)(m-1)
 \longrightarrow 0
\end{equation}
\begin{equation}
 0\longrightarrow \left(g^*\Omega^1_{\P^n}\right)(m-1) \longrightarrow
 \O_\Sigma(-1)^{\oplus n+1} \longrightarrow \O_\Sigma(m-1) \longrightarrow 0
\end{equation}
Since $H^0(\Sigma, \O_\Sigma(-1))=0$, (2.13) yields
$H^0\left(\Sigma, \left(g^*\Omega^1_{\P^n}\right)(m-1)\right)=0$.
Therefore (2.12) implies $H^0\left(\Sigma,
\left(g^*\Omega_Y^1\right)(m-1)\right)=0$ because $H^1(\Sigma,
\O_\Sigma(m(1-e)-1))=0$.  Hence $n-\delta +m(e-n)-1<0$ by (2.9),
as desired.
\end{proof}
\end{prop}

\subsection*{Proof of Theorem 2.}
Suppose that $F^{-1}(Y)\neq X$.  Then
$H$ is not empty, and
Proposition 2.4 implies
$$
n+m(e-n) \leq \deg \Sigma \leq \deg H = em-d.
$$
Therefore $d\leq n(m-1)$.  If $Y$ is of general type, i.e. $e\geq n+2$, then Proposition 2.2(4)
implies $d>mn$.  This contradiction finishes the proof.

\subsection*{Remark}
Suppose the ground field $k$ has positive characteristic.  If the characteristic is large, say $\text{char}\hspace{0.04 cm} k> \alpha$, where
$\alpha:= \frac{em-d}{e}m^{n-2}$, then the morphism
$F|_\Sigma:\Sigma\rightarrow F(\Sigma)$ is separable.  Indeed, the Grothendieck-Lefschetz Theorem, \cite[Theorem 4.3.2]{H1}, implies that the divisor $F(\Sigma)\subset Y$ is the intersection of $Y$ with another hypersurface.  So
one can check that  $\deg
F|_\Sigma\leq \alpha$ by applying the projection
formula to $F|_\Sigma\rightarrow F(\Sigma)$.

It follows that if $\text{char}\hspace{0.04
cm} k> \alpha$, then the proof of Proposition 2.4 is still valid.  Hence,
Theorem 2 will also hold in positive characteristic if char\hspace{.04cm}$k >\alpha$.

\begin{cor}
If $f:X\rightarrow Y$ is a nonconstant morphism between hypersurfaces in $\P^n$, $n\geq 4$, such that $Y$ is nonsingular and of general type, then $\deg Y$ divides $\deg X$ with quotient $m$ such that $f^*\O_Y(1)=\O_X(m)$.
\end{cor}
\begin{proof}
By Theorem 2, there is a morphism $F:\mathbb P^n\rightarrow \mathbb P^n$ such that $X=F^{-1}(Y)$ and $F^*\O(1)=\O(m)$.  It follows that $X$ is a hypersurface of degree $m\cdot\deg Y$.
\end{proof}

\subsection{Hypersurfaces Not of General Type}
The ground field will have characteristic zero unless indicated otherwise.
We will show that if $3\leq e\leq n+1$ and $d$ is not too much larger than $e$, then the conclusion of Theorem 2 still holds. The following definition will be central to our point of view.

\subsection*{Definition}
If $Z$ is any scheme and $F:Z\dashrightarrow \P^n$ is a rational map given by
sections $F_0, \dots, F_n$ of some line bundle on $Z$, then let
\emph{indet$(F)$} denote the scheme of common vanishing of the $F_i$ in $Z$:
$$indet(F):= V(F_0,\dots, F_n)\subset Z.$$

\begin{lem}
Using the notation at the beginning of Section 2: \\
(1) If $H\neq 0$, then $\deg H=em-d\geq e$.  This holds for $e\geq 2$.\\
(2) If $p\in\indet(F)$, then $H$ has order at least $e$ at $p$,
regardless of the characteristic of the ground field.
\begin{proof}
Suppose $p\in\indet(F)$ is a reduced closed point. If $Y=V(G)$ for
a homogeneous polynomial $G=G(y_0,\dots,y_n)$ of degree $e$, then
$F^{-1}(Y):=V(G(F_0, \dots, F_n))$ has order at least $e$ at $p$
because the $F_i$ are all zero at $p$. But $F^{-1}(Y)=X+H$, and
$p$ is not contained in $X$. So $H$ has order at least $e$ at $p$.
This proves (2), and it proves (1) in case $F$ is not defined at
some point of $H$.

If $F|_H$ is a morphism, then (1) follows from Proposition 2.2(2) in
case $e\geq 3$.  And in the case $e=2$, we need only see that
$d\neq 1$.  However, Lazarsfeld shows in \cite{L} that if a smooth
variety $Y$ is the image of a morphism from a projective space,
then $Y$ is itself a projective space.
\end{proof}
\end{lem}

\begin{prop}
If $m=1,2$, then the conclusion of Theorem 2 holds, i.e. there is a morphism $F:\P^n\rightarrow\P^n$ such that $f=F|_X$ and $X=F^{-1}(Y)$.
\begin{proof}
If $m=1$, then the image of $F:\Pn\dashrightarrow\Pn$ is a linear
subspace that contains $Y$.  So $F$ is an automorphism of $\Pn$,
and $d=e$.

Suppose $m=2$ and $d\neq 2e$.  Then $X$ and $H$ both have degree
$e$ by Lemma 2.6 and Proposition 2.2(2).  If $e\geq 3$, then
$m=1$ by Proposition 2.2(3), which is a contradiction.
If
$e=2$, then $m=1$ because every nonconstant morphism of smooth
quadrics in $\P^n$ is an isomorphism for $n\geq 4$, cf. \cite{PS}.
This contradiction shows $d=2e$ after all.
\end{proof}
\end{prop}

\begin{prop}
Fix $d$, $e$, $m$ with $e\geq 3$, and assume one of the following three conditions holds: \\
(i)  $d<e^2$\\
(ii) $d>(m-1)^2$ \\
(iii) $m\leq e$ \\
Then the conclusion of Theorem 2 holds for $n$ sufficiently large.
\begin{proof}
Using Theorem 1 and formulas (2.3) and (2.4), let $n$ tend
to infinity and get $d-1\geq m(e-1)$.  If $H\neq 0$,
then $em-d\geq e$ by Lemma 2.6(1).  Together, these two
inequalities contradict each of the three conditions above.
\end{proof}
\end{prop}

\subsection*{Examples in Characteristic Zero.}
Theorem 1 gives an upper bound on the polynomial degree $m$ of $f$ whenever $e\geq 3$.  Using these explicit upper bounds, along with Proposition 2.4 and Lemma 2.6(1), one can check that
the conclusion of Theorem 2 holds for the following cases in $\P^4$:
\begin{align*}
 e=3 \ \ \ \ \ & d\leq 4  \\
 e=4 \ \ \ \ \ & d\leq 10                    \\
 e=5 \ \ \ \ \ & d= 1,\dots, 23, 25, 26, 29
\end{align*}

\subsection*{Examples in Positive Characteristic.}  Theorem 1 and formulas (2.3), (2.4) hold in arbitrary characteristic.  So we may compute upper bounds on the polynomial degree $m$ of $f$ in positive characteristic as well.

If $em\neq d$ and $\deg H=1$, then $F|_H$ is a morphism by Lemma 2.6(2).  This  is impossible
because $Y$ is not the image of a morphism from $\P^{n-1}$, as
shown in \cite{L}.  So if $H\neq 0$, then $\deg H>1$.

Using the fact $em-d>1$ and the explicit upper bounds on $m$ that we obtain from Theorem 1, we see
that the conclusion of Theorem 2 holds for the following cases in $\P^4$:
\begin{align*}
 e=3 \ \ \ \ \ & d\leq 3                    \\
 e=4 \ \ \ \ \ & d\leq 8                    \\
 e=5 \ \ \ \ \ & d = 1,\dots, 11, 14                   \\
 e=6 \ \ \ \ \ & d = 1,\dots, 14, 17, 18             \\
 e=7 \ \ \ \ \ & d= 1,\dots, 17, 20, 21, 22, 27
\end{align*}

\subsection*{Question} One can ask if the general type hypothesis of Theorem 2 is too strong.  The results of Section 2.3 seem to indicate that this is indeed the case.  To be precise, if $f:X\rightarrow Y$ is a nonconstant separable morphism of hypersurfaces in $\P^n$, $n\geq 4$, such that $Y$ is nonsingular of degree at least 2, is it true that there is an endomorphism $F:\P^n\rightarrow\P^n$ such that $f=F|_X$ and $X=F^{-1}(Y)$?

\bibliographystyle{amsplain}
\bibliography{xbib}

\begin{thebibliography}{}

\fontsize{10}{12}\selectfont

\bibitem{A}  Amerik, E. :  On a Problem of Noether-Lefschetz type,
Comp. Math. 112 (1998), 255-271
\bibitem{ARV}  Amerik, E., Rovinsky, M., Van de Ven, A. :  A boundedness theorem for
morphisms between threefolds, Ann. Inst. Fourier 49 (1999),
405-415
\bibitem{B} Beauville, A. :  Endomorphisms of hypersurfaces
and other manifolds, preprint: arXiv (2000)
\bibitem{CG} Clemens, C. H., Griffiths, P. A. : The intermediate
Jacobian of the Cubic Threefold, Annals of Mathematics, 95 (1972),
281-356
\bibitem{E} Eisenbud, D. :
Commutative Algebra with a View Towards Algebraic Geometry, Grad.
Texts in Math 150, Springer-Verlag (1995)
\bibitem{F} Fulton, W. : Intersection Theory,
$2^{nd}$ Ed., Springer-Verlag (1998)
\bibitem{HRS}  Harris, J., Roth, M.,
Starr, J. :  Curves of Small Degree on Cubic Threefolds, preprint:
arXiv (2001)
\bibitem {H1} Hartshorne, R. :  Ample Subvarieties of Algebraic
Varieties, Lecture Notes in Math. 156, Springer-Verlag (1970)
\bibitem {H2} Hartshorne, R. :  Algebraic Geometry, Graduate Texts in Math 52, Springer-Verlag
(1977)
\bibitem {L} Lazarsfeld, R. :  Some applications of the theory of
Positive vector bundles, Lecture Notes in Math. 1092, (1984)
\bibitem {PS} Paranjape, Srinivas :  Self maps of Homogeneous
Spaces, Inven. Math. 98 (1989), 425-444
\bibitem {S1} Schuhmann, C. : Mapping threefolds onto three dimensional quadrics,
Mathematische Annalen 306 (1996), 571-581
\bibitem {S2} Schuhmann,
C. : Morphisms between Fano Threefolds, J. Algebraic Geometry 8
(1999), 221-244
\end{thebibliography}
\end{document}